\documentclass[ amscd, syntonly,  psfig,  amssymb, 11pt]{amsart}
\usepackage{psfig}
\newtheorem{theo}{Theorem}[section]
\newtheorem{pro}[theo]{Proposition}
\newtheorem{lem}[theo]{Lemma}
\newtheorem{defi}[theo]{Definition}
\newtheorem{remark}[theo]{Remark}
\newtheorem{corol}[theo]{Corollary}

\def \proof {{\bf Proof$\colon$}\ }

\def \Z{{\bf Z}}
\def \Q{{\bf Q}}
\def \R{{\bf R}}
\def \N{{\bf N}}
\def \ad {{$K\buildrel n\over \longrightarrow K'$}}

\def\no{\noindent}
\begin{document}

\title[]{Knot adjacency and fibering}

\author[E. Kalfagianni]{Efstratia Kalfagianni$^1$}
\author[X.-S. Lin]{Xiao-Song Lin$^2$}

\address[]{Department of Mathematics , Michigan State
University,
E. Lansing, MI, 48824}

\email[]{kalfagia@math.msu.edu }

\thanks{{\it AMS classification numbers:} 57M25, 57M27, 57M50.}


\thanks{type invariants, symplectic structures.}
\thanks{$^1,  ^2$ The research of the authors is partially supported by the NSF}

\address[]{Department of Mathematics , 
University of California,
Riverside, CA, 92521}
\email[]{xl@math.ucr.edu}

\begin{abstract} It is known that the
Alexander polynomial detects fibered knots and 3-manifolds
that fiber over the circle. In this note, we show that when the
Alexander polynomial becomes inconclusive, the notion
of {\em knot adjacency}
can be used to obtain obstructions
to fibering of knots and of 3-manifolds. As an application, given a
fibered knot $K'$, we construct infinitely many non-fibered knots
that share the same Alexander module with $K'$. Our construction
also provides, for every $n\in N$, examples of
irreducible 3-manifolds that cannot be distinguished
by the Cochran-Melvin finite type invariants of order $\leq n$.

\end{abstract}

%\maketable
\smallskip
\maketitle
\tableofcontents{}

\medskip

\section{Introduction} The problem of detecting fiberedness (or non-fiberedness) of knots,
has been studied considerably from both the algebraic
and the geometric topology viewpoint.
The classical Alexander polynomial of a fibered
knot is known to be {\em monic} and this provides an effective criterion for detecting
fibered knots.
The converse is not in general true, although it holds
for several special classes of knots including alternating knots and knots
up to ten crossings. Non-commutative generalizations of the Alexander
polynomial, such as the higher order Alexander polynomials defined
in \cite{kn:co} and suitable versions of the
twisted Alexander polynomials are known to
detect non-fibered knots with monic Alexander polynomial
(\cite{kn:gkm, kn:cha, kn:fk}). A geometric procedure to
detect fibered
knots was developed by Gabai
in \cite{kn:ga2}.

The purpose of this paper is to introduce a new criterion for
detecting non-fibered knots when the Alexander polynomial fails,
and to present several applications of this criterion. Our approach
combines both the algebraic and the geometric point of view.
To describe
our results, we recall that knot $K$ is called $n$-adjacent to another
knot $K'$, if $K$
admits a projection containing $n$ crossings
such that changing any $0<m\leq n$ of them yields a
projection of $K'$.
The notion of knot adjacency was studied in \cite{kn:kl1} and the theory was further developed 
in \cite{kn:kl}, \cite{kn:k1}. In particular, in \cite{kn:k1} we
showed that high degree adjacency ($n>1$) to fibered knots
imposes strong restrictions on the knot genus. In this paper,
we explore the role of knot adjacency as obstruction to fiberdness. We show
that when the Alexander polynomial provides
inconclusive evidence, high degree knot adjacencies obstruct
knots to be fibered. More precisely, we have the following:

\begin{theo} \label{theo:detect1} Let $K, K'$ be distinct knots
with equal Alexander polynomials.
Suppose that $K'$ is fibered.
If $K$ is $n$-adjacent to $K'$ for some $n>1$,
then $K$ is not fibered. Furthermore,
the 3-manifold $K(0)$
obtained by 0-Dehn surgery of $S^3$ along $K$
does not admit  a fibration over $S^1$.
\end{theo}

Since it is known that all  knots up to ten crossings 
that have monic Alexander polynomial
are fibered (\cite{kn:ga2}), the knots that are detected by
Theorem \ref{theo:detect1}
have large crossing numbers.
In Section 2
we give a systematic procedure that constructs vast
families of such knots.
In particular,
given a
fibered knot $K'$, we construct infinitely many non-fibered knots
that share the same Alexander module with $K'$.
To state our results,
for a knot $K$, let $\Delta_K(t)$ denote the Alexander
polynomial of $K$,
let $\delta(K)$ denote the degree of $\Delta_K(t)$
and let $g(K)$ denote the genus of $K$.  
It is well known that if
a knot $K$ is fibered, then $m(K):= 2g(K)-\delta(K)=0$.  
\vfill
\eject

\begin{theo} \label{theo:fibered0} Given a fibered knot $K'$
there exist infinitely many knots
$\{K_n\}_{n\in \N}$, such that:

\no (a) \ $K_n$ and $K'$
have the same Alexander module.

\no (b)\ $g(K_n)> g(K')$ and ${\displaystyle {{{n+3}\over 6}}\leq g(K_n)}$.

\no (c) \ $m(K_n)>0$ and $\lim_{n \to \infty}m(K_n)=\infty$. In particular, $K_n$ is not fibered.

\no (d)\  $K_n(0)$ does not fiber over $S^1$.

Furthermore, if $K'$ is a prime knot then $K_n$ can be chosen to be prime.
\end{theo}

In the course of the proof of Theorem \ref{theo:fibered0},
we show that all the knots $\{K_n\}_{n\in \N}$ share a common
equivalence class of Seifert matrices with $K'$.
Thus,
$K_n$ and $K$ cannot be distinguished by {\em abelian}
invariants (e.g.  torsion numbers, signature,
Blanchfield linking forms). This gives a new proof
of
a slightly weaker
version of the main result of \cite{kn:cha}.
\smallskip

As another application of Theorem \ref{theo:detect1},
we construct families of
3-manifolds that cannot be distinguished by certain
finite type invariants. For ${\Z}$-homology 3-spheres
these invariants were defined by Ohtsuki in \cite{kn:o}.
An extension of Ohtsuki's theory to arbitrary
3-manifolds was proposed by Cochran and Melvin in \cite{kn:cm}.
There exist constructions that yield  irreducible
distinct homology 3-spheres with the same finite type invariants
of bounded order (\cite{kn:k}).
Here we give examples of  irreducible
3-manifolds with non-trivial homology that cannot be distinguished
by their finite type invariants in the sense of \cite{kn:cm}.
To state our results, for a knot $K$ and a number $s\in \Q$,
let $K(s)$ denote the 3-manifold obtained by $s$-Dehn surgery of $S^3$ along $K$.
We show that if $K,K'$ are knots such that $K$ is $n$-adjacent to $K'$, then, for every $s\in \Q$,
$K(s)$ and $K'(s)$ cannot be distinguished
by any finite type invariant of order $<n$ (Proposition \ref{pro:same}). Combining this with
Theorem \ref{theo:detect1}
we obtain the following:

\begin{corol} \label{corol:S2xS1}
For every $n\in \N$, there exist closed, irreducible 3-manifolds
$M, M'$ such that:

\no (a)\ $M$ and $M'$ have the homology type of $S^2\times S^1$.

\no (b) \ For every commutative ring $\mathcal R$ with unit,
$M$ and $M'$
have the same
$\mathcal R$-valued Cochran-Melvin finite type invariants of
order $< n$.

\no (c) \ Exactly one of $M, M'$ fibers over $S^1$.

\end{corol}

The examples of Corollary \ref{corol:S2xS1}
are obtained by 0-Dehn surgery on knots $K\subset S^3$.
If $0\neq s:={\displaystyle {a\over b}}$, then
$K(s)$ is a $\Q$-homology sphere. In particular,
if $K$ is the unknot then $K(s)$ is the Lens space
$L(a, \ b)$ and thus $\pi_1(K(s))=\Z _a$.
Combining our work with the Cyclic Surgery Theorem
of Culler, Gordon, Luecke and Shalen (\cite{kn:cgls}) we obtain the following:

\begin{corol} \label{corol:lens}
Given $n\in\N$, there exists a knot $K\subset S^3$ and an
integer $r_K$ such that for every $s:={\displaystyle {a\over b}}\neq r_K, r_K+1$ we have:

\no (a)\ $\pi_1(K(s))\neq \pi_1(L(a, \ b))$.

\no (b) For every commutative ring $\mathcal R$ with unit,
$K(s)$ and $L(a, \ b)$ have the same
$\mathcal R$-valued Cochran-Melvin finite type invariants of
order $< n$.
\end{corol}

The examples of Corollary \ref{corol:S2xS1} and Corollary \ref{corol:lens}
are the first known examples of irreducible,
non $\Z$-homology spheres
that are not distinguished by
finite type invariants of bounded order.

The paper is organized as follows: In Sections 2 and 3 we state a fibering
criterion in terms of knot adjacency and discuss its applications to detecting non-fibered knots and 3-manifolds.
In Section 4 we begin by recalling
from \cite{kn:cm}
the definition of finite type invariants
for arbitrary 3-manifolds. Then, we apply our results from
Sections 2 and 3 to construct examples of 3-manifolds
with the same finite type invariants of bounded orders
and to prove Corollaries \ref{corol:S2xS1}
and \ref{corol:lens}. In Section 5, we give the proof of Theorem 1.3.
Section 6  is a short appendix in which we discuss that,
combined with work of Kronheimer,
Theorem
\ref{theo:fibered0} can be used
to detect the existence of symplectic structures on certain
4-manifolds.
\vskip 0.04in
\smallskip

{\bf Acknowledgment.} Effie Kalfagianni  thanks the Institute for Advanced Study
for their hospitality while part of the research described in this paper was
completed and for partial support through a
research grant.
She also acknowledges the partial support of the NSF through
grants DMS-0306995  and FRG/DMS-0456155.
Xiao-Song Lin acknowledges the partial support of the NSF through
grants DMS-0404511 and FRG/DMS-0456217.

\medskip

\section{Adjacency to fibered knots and the Alexander polynomial}

Let $K$ be a knot in ${\bf S}^3$ and let $q\in \Z$. 
A generalized crossing of order $q$
on a projection of $K$ is a set $C$
of $|q|$ twist crossings
on two strings that inherit
opposite orientations from any orientation of $K$.
If $K'$ is obtained from $K$ by changing
all the crossings in $C$ simultaneously, we will say that
$K'$ is obtained from $K$ by a {\em generalized crossing change} of order $q$. 
Note that if $|q|=1$, $K$ and $K_1$ differ by an ordinary crossing change
while if $q=0$ we have $K=K'$.
A {\it crossing disc} corresponding to a generalized crossing $C$ of a
knot $K$ 
is an embedded disc $D\subset S^3$
such
that $K$ intersects ${\rm int}(D)$ twice, with
zero algebraic intersection number (once for each string of $K$ forming the crossing). The curve $\partial D$ is called a crossing
circle corresponding to $C$.
The crossing is called {\em nugatory} if $\partial D$ 
bounds disc in the complement of $K$.

\begin{defi} \label{defi:adj}
We
will say that $K$ is $n$-adjacent to $K'$, for some
$n>0$,
if $K$ admits a projection containing $n$ generalized crossings
such that changing any $0<m\leq n$ of them yields a
projection of $K'$. We will write \ad.
\end{defi}

For a knot $K$, let
$g(K)$ denote its genus  and let $\delta:=\delta(K)$
denote the degree of its
Alexander polynomial $\Delta_K(t)=\sum_{i=0}^{\delta} a_i t^i$.
We have the following theorem, that contains the first conclusion in 
Theorem \ref{theo:detect1}.

\begin{theo}\label{theo:detect}
Let $K,K'$ be knots such that $\delta(K)=\delta(K')$ and
\ad , for some $n>1$.
If $K'$ is fibered, then either $K$ is isotopic to $K'$ or
we have  $g(K)>g(K')$. Furthermore, in the later case,
$K$ is neither fibered nor alternating.
\end{theo}
\proof Since $K'$ is fibered, we have $\delta(K')=2g(K')$. Thus
we have 
$$\delta(K)=\delta(K')=2g(K'). \eqno (1)$$
The assumption that \ad \ allows as to apply the results of \cite{kn:kl1, kn:k1}:
Indeed, as it was shown in \cite{kn:k1}, if $K$ is $n$-adjacent to a
fibered knot $K'$, for some $n>1$,
then
either $K$ is isotopic to $K'$ or
we have $g(K)>g(K')$. If $g(K)>g(K')$,
then (1) implies that $m(K):=2g(K)-\delta(K)>0$
and thus $K$ is not fibered or alternating. \qed
\smallskip

By Corollary 8.19 of \cite{kn:ga}, $K(0)$ fibers over
$S^1$ precisely when $K$ is a fibered knot. Thus, the second conclusion of 
Theorem \ref{theo:detect1} holds.
The next theorem and its proof provide vast collections of non-fibered knots with monic
Alexander polynomials that are detected by
Theorem \ref{theo:detect}.

\begin{theo} \label{theo:construct} Given a fibered knot
$K'$ and $n>1$,  there exists a
knot $K$ with the following
properties:

\no (a)\  \ad.

\no (b)\  $K$ and $K'$
have the same Alexander module.

\no (c)\  $g(K)>g(K')$.

\no (d)\ If $K'$ is a prime knot, then $K$ can be chosen to be prime.
\end{theo}

\begin{remark} \label{n=1}{\rm  Theorem \ref{theo:construct} (a)-(c) remains
true 
for $n=1$. To see this, let $J$ be a knot with trivial Alexander polynomial 
that can be unknotted by a single generalized crossing change (e.g. an
untwisted Whitehead double of any knot). Then the connected sum $K:=
J\#K'$ has the properties (a)-(c) of Theorem \ref{theo:construct}.}
\end{remark}
\smallskip

The proof of Theorem \ref{theo:construct} and that of
Theorem \ref{theo:fibered0} are given in section 5. In the next two sections we present
applications of Theorems \ref{theo:detect1} and \ref{theo:fibered0}.
\medskip

\section{Obstructing fibrations}

Let ${\mathcal K}_1$ denote the set of isotopy classes of knots with monic Alexander polynomial  and recall that  for $K\in {\mathcal K}_1$,  we denote $m(K):= 2g(K)-\delta(K)\geq 0$, where
$\delta(K)$ is the degree of the  Alexander polynomial of $K$.
For  $K\in {\mathcal K}_1$, let $F_K \subset {\mathcal K}_1$ denote the set of
isotopy classes of  fibered knots, that are distinct from $K$
but share the same Alexander polynomial with $K$. As shown in (\cite{kn:morton},
if $\delta(K)>2$ then $F_K$ is infinite. On the other hand, it is known that the only fibered knots of genus $\leq 1$
are the trefoils, the figure eight and the unknot. It is easy to see that
all degree $\leq 2$ monic Alexander polynomials are realized by these knots. Thus,
if $\delta(K)\leq 2$, then $F_K$ is finite. Fix $K\in {\mathcal K}_1$.
For $K'\in F_K$, let 
$$I(K, K'):= \{ n>1 \ \ | \ \ K\buildrel n\over \longrightarrow K' \}.$$
We define
$a(K, \ K'):= {\rm max} \{ n\in I(K, K')\},$
if the set $I(K, K')$ is non-empty. Otherwise,
we define 
$a(K, \ K'):=0$.
Finally, we define $$\alpha(K):= {\rm max} \{ a(K, \ K'), \  K'\in F_K \}.$$
The quantity $\alpha(K)$ is clearly an invariant of $K$. Roughly speaking, it  measures the degree of adjacency of a knot with monic Alexander polynomial
to fibered knots with the same polynomial. By definition, we have $\alpha(K)\in {\bf N} \cup \{\infty\}$; the following proposition shows that, in fact,
$\alpha(K)< \infty$, for every $K\in {\mathcal K}_1$.

\begin{pro} \label{pro:invariant} The invariant $\alpha$
has the following properties:

\no (a)\ We have $0\leq \alpha(K)\leq 6g(K)-3$,
for every $K\in {\mathcal K}_1$.

\no (b)\ If $m(K)=0$ then $\alpha(K)=0$. In particular,
if
$K$ is fibered then $\alpha(K)=0$.
\end{pro}
\proof Clearly we have  $\alpha(K)\geq 0$.
By definition, if $\alpha(K)>0$
then there is a fibered knot $K'\neq K$ such that
$\Delta_K(t)=\Delta_{K'}(t)$ and \ad , for some $n>1$.
By Theorem \ref{theo:detect}, we must have
$g(K)>g(K')$ and Theorem 1.3 of \cite{kn:kl1} applies to conclude that
$n\leq 6g(K)-3$.  Thus part (a) is proved.
To see part (b) suppose that $m(K):=2g(K)-\delta(K)=0$ and
that $\alpha(K)>0$.
Then, by definition, there is a fibered knot $K'\neq K$
such that $\delta(K')=\delta(K)$ and \ad , for some $n>1$.
Since $K'$ is fibered we have $\delta(K')=2g(K')$. But since
$\delta(K)=2g(K)$,
we conclude that $g(K)=g(K')$.
But this is impossible since by Theorem \ref{theo:detect}, we must have
$g(K)>g(K')$. \qed
\smallskip

The proof of the next  corollary uses
$\alpha(K)$  to produce
infinitely many non-fibered knots 
with a given monic Alexander polynomial.
\begin{corol} \label{corol:infinitely}
For every fibered knot $K'$ there exist infinitely many non-fibered knots $\{ K_n\}_{n\in N}$
each of which has the same Alexander module
with $K'$. Furthermore, if $K'$ is a prime knot then $K_n$ can be taken to be prime.
\end{corol}
\proof Let $K'$ be a fibered knot and fix $n'>1$.
By Theorem \ref{theo:construct}
there  exists a  knot $K_1$
such that $$\Delta_{K_1}(t)=\Delta_{K}(t), \ \
K_1\buildrel n'\over \longrightarrow K', \ \ {\rm and}\ \ g(K_1)>g(K').$$
It follows that
$2g(K_1)>\delta(K_1)$ which implies that
$K_1$ is non-fibered.
Clearly $\alpha(K_1)\geq \alpha(K_1, \ K')\geq n'$.
Suppose, inductively, that we have constructed
non-fibered knots $K_1, \ldots, K_m$ such that
$\alpha(K_m)>\ldots >\alpha(K_1)\geq n'$ and
$\Delta_{K_m}(t)=\ldots =\Delta_{K_1}(t)=\Delta_{K}(t)$.
Clearly, $K_1, \ldots, K_m$ are distinct.
Now choose $n>> \alpha(K_m)$ and let $K_{m+1}$
be any knot obtained by applying Theorem \ref{theo:construct} to this $n$. \qed
\smallskip

\smallskip

By Theorem \ref{theo:fibered0}(b), the knots
$\{K_n\}_{n\in \N}$ can be chosen so that
$g(K_{n+1})> g(K_{n})$. By \cite{kn:ga}, the 3-manifold $K_n(0)$ contains
a closed, embedded, orientable, non-separating surface of genus
$g(K_n)$ and contains no such surface of smaller genus.
It follows that the manifolds $\{K_n(0)\}_{n\in \N}$
are all distinct. On the other hand, since the Alexander module 
of $K_n(0)$ is the same as that of $K_n$, all these 3-manifolds have the same Milnor torsion
(\cite{kn:tu}). Thus we obtain the following:

\begin{corol}\label{corol:milnor}
Given a fibered knot $K'$
there exist infinitely many non-fibered knots
$\{K_n\}_{n\in \N}$, such that the 3-manifolds
$\{K_n(0)\}_{n\in \N}$ are all distinct but have the same Milnor torsion
with $K'(0)$.
\end{corol}

\medskip

\section{Cochran-Melvin invariants of 3-manifolds}
There exist several constructions of
$\Z$-homology 3-spheres that are indistinguishable
by Ohtsuki invariants of bounded order.
The examples given in
\cite{kn:k} are obtained by Dehn surgery along suitable knots in $S^3$.
In \cite{kn:cm} Cochran  and Melvin generalized Ohtsuki's
theory to define finite type invariants for arbitrary
3-manifolds.
The knots we construct
in this paper, and in particular these in the proof of Theorem
\ref{theo:construct},
fit nicely into
the theory of \cite{kn:cm} and lead to
a natural extension of the construction
of \cite{kn:k} in this setting.
Before we state our results we recall some
definitions.

\begin{defi} \label{defi:admissible }A  framed link $L $
in a closed, oriented 3-manifold $N$
is called {\em admissible} iff we have:

 \no\ (i) Each component of $L$ is null-homologous
in $N$.

\no \ (ii) All the pairwise linking numbers of $L$ in $N$ vanish.

\no \ (iii) The framings are $\pm 1$ with respect to
the longitudes given by (i).
\end{defi}

Let $N$ be a closed oriented 3-manifold.
The set ${\mathcal S}:={{\mathcal S}}( N)$ of
homeomorphism
classes of 3-manifolds that are $H_1$-cobordant to $N$ is precisely
the set of 3-manifolds obtained by surgery of $N$ along admissible links
(\cite{kn:cm}). Let ${\mathcal R}$ be a commutative ring with
unit, and let ${\mathcal M}(N)$ 
be the ${\mathcal R}$-module freely
spanned by $\mathcal S$.
For $M\in {\mathcal S}$ and an admissible link $L\subset M$  
define $[M, \ L] \in {\mathcal M}(N)$
by
$$[M, \ L]:= \sum_{L^{'} \subset L} (-1)^{\# L^{'}} M_{L^{'}}\eqno(6)$$
where $L^{'}$ ranges over all sublinks of $L$
(including the empty one). Here
$\# L^{'}$ denotes the number of components of $L^{'}$
and $M_{L^{'}}$ denotes the 3-manifold
obtained from $M$ by surgery along $L^{'}$.
For $l\geq 0$, let  ${{\mathcal M}_{l}}(N)$ denote the submodule of ${\mathcal M}(N)$
that is freely spanned by all expressions $[M, \ L]$,
where $M\in {\mathcal S}$ and $L$ is an  admissible link in $M$
with $\# L \geq l$. Let $\mathcal H$ denote the set of
$H_1$-cobordism classes of closed, oriented 3-manifolds;
for $i\in {\mathcal H}$ choose a representative $N_i$.
Let

$${\mathcal M}:=\bigoplus_{i\in {\mathcal H}} {\mathcal M}(N_i)\ \ {\rm and} \ \
{{\mathcal M}_{l}}:= \bigoplus_{i\in {\mathcal H}} {\mathcal M}_{l}(N_i).$$
\smallskip

\begin{defi} \label{defi:cocm} \ {\rm (\cite{kn:cm}) }\ A
functional $f: {\mathcal M}/ {{\mathcal M}_{n+1}}
\longrightarrow \mathcal R$ is
called an $\mathcal R$-valued
finite type
invariant
of order $\leq n$. We will use ${\mathcal F}_n$ to denote
the space of all such functionals.
\end{defi}

One can see that
$${\mathcal F}_n\cong \bigoplus_{i\in {\mathcal H}}{\rm Hom}({\mathcal G}_{n}(N_i)
, \ {\mathcal R})\ \  
{\rm where} \ \ {\mathcal G}_{n}(N_i):={\mathcal M}(N_i)/ {{\mathcal M}_{n+1}(N_i)}.$$
Thus, the invariants of finite type of
 \cite{kn:cm} are constructed from invariants in each $H_1$-cobordism class.
Moreover, the invariants of type $0$ are exactly the functionals
${\mathcal H}\longrightarrow {\mathcal R}$.
In \cite{kn:cm} it is shown that, for every $n\in\N$,
${\rm Hom}({\mathcal G}_{n}(N_i)
, \ {\mathcal R})$ is a finite dimensional
non-trivial $\mathcal R$-module. To state our results, for a
knot $K\subset S^3$ and a rational number $s\in \Q$, let
$K(s)$ denote the 3-manifold obtained by $s$-surgery of $S^3$ along $K$. Note that $K(s)$ is either a rational homology 3-sphere or a homology $S^2\times S^1$ manifold.

\begin{pro} \label{pro:same}
Suppose that $K,K'$ are knots such that
\ad , for some $n>0$. Suppose, moreover, that
there exists a collection of $n$ ordinary crossings
that exhibit $K$ as $n$-adjacent to $K'$.
Then, 
for every $s\in \Q$, we have:
$$f(K(s))=f(K'(s)),$$
for every $f\in {\mathcal F}_{n-1}$.
\end{pro}
\proof 
Fix $n>0$ and  let $K$, $K'$  be knots
such that $K$ admits a collection 
of ordinary crossings $\mathcal C$
that exhibit it as $n$-adjacent to $K'$.
Let $L\subset S^3$ be an $n$-component  link consisting of
a crossing circle for each of the crossings in $\mathcal C$. The crossing change
can be achieved by doing surgery of $S^3$ along the corresponding
crossing circle; the framing of the surgery is $+1$
or $-1$ according to whether the crossing is positive or negative. Thus $L$ can be considered as admissible.
Since the linking number of $K$ with each component of $L$ is zero,
each component of $L$ is null-homologous in $S^3\setminus K$.
Since, for every $s\in \Q$, there is an epimorphism
$H_1(S^3\setminus K)\longrightarrow H_1(K(s)),$
it follows that the image of $L$
in $M:=K(s)$ is an admissible link; we will still denote this link by $L$. Let $L'\subset L\subset M$ be a
non-empty sublink of $L$ in $M$. The 3-manifold $M_{L'}$,
obtained by surgery of $M$ along $L'$, can be alternatively
described as follows: First perform surgery of $S^3$ along
$L'$; this gives back $S^3$ but it changes $K$ to $K'$.
Then, perform $s$-Dehn surgery of $S^3$ along $K'$.
From these considerations we conclude that
$M_{L'}=K'(s)$. Thus (6) yields
$$[M, \ L]=M-K'(s).\eqno(7)$$
Now let $f\in {\mathcal F}_{n-1}$. Since by definition
$f([M, \ L])=0$, from (7) we obtain $f(M)=f(K'(s))$ as desired. \qed
\smallskip

\begin{corol} \label{corol:samef}
Let $n>0$. For every fibered knot $K'\subset S^3$ there exists a
non-fibered knot $K\subset S^3$ such that 
for every $s\in \Q$, we have:
 $$f(K(s))=f(K'(s)),$$
for every $f\in {\mathcal F}_{n-1}$.
\end{corol}

\proof Given $n$ and $K'$ as above,
let $K$ be a knot corresponding to $n$ and $K'$ in the sense of Theorem
\ref{theo:construct} if $n>1$. Since we do not require $K$ to be prime,
we will use Remark \ref{n=1} to conclude that such a $K$ also exists when
$n=1$. By the same token, the proof of Theorem
\ref{theo:construct} shows that we can choose $K$ so that
it is shown to be $n$-adjacent to
$K'$ by a collection of ordinary crossings.
Thus the corollary follows from Proposition \ref{pro:same}. \qed
\smallskip
\smallskip

By Gabai's work (\cite{kn:ga}), $K(0)$ is irreducible if
$K$ is
non-trivial and it fibers over $S^1$ precisely when $K$ is fibered.
Combining these facts with  Corollary \ref{corol:samef} and Theorem \ref{theo:fibered0}
we obtain Corollary \ref{corol:S2xS1}.
Next we give the proof of Corollary \ref{corol:lens}.
\smallskip
\smallskip

{\bf Proof of Corollary \ref{corol:lens}:} By
Theorem
\ref{theo:construct}, for every
$n\in \N$, there exists a
non-trivial knot $K$ that is $n$-adjacent to the trivial knot
and has trivial Alexander polynomial.
Part (b) of the corollary follows immediately from 
Proposition \ref{pro:same}.  Corollary 1
of 
\cite{kn:cgls} states that if $K$ is not a
torus knot then only for integer slopes $r$,
$K(r)$ can have cyclic
fundamental group. Furthermore, there can be at
most two such integers and if there are two they have to be successive.
Since $K$ has trivial Alexander polynomial
it cannot be a torus knot. Thus 
part (a) follows immediately from
Corollary 1 of 
\cite{kn:cgls}. \qed

\medskip

\section{Constructing the knots $K_n$}
In this section we give the proof of Theorem \ref{theo:construct}.
First let us explain how
Theorem \ref{theo:fibered0} follows from Theorem \ref{theo:construct}:
Given a fibered knot $K'$, for every $n\in \N$,
let $K_n$ be a non-fibered knot guaranteed by Theorem \ref{theo:construct}.
Since $K_n$ is $n$-adjacent to $K'$, by \cite{kn:kl1},
we have $n\leq 6g(K_n)-3$. Hence ${\displaystyle {{{n+3}\over 6}}\leq g(K_n)}$
and parts (a)-(c) of Theorem \ref{theo:fibered0} follow.
For part (d),
we repeat that it follows by Corollary 8.19 of \cite{kn:ga}.
\smallskip

 Before we can proceed with the proof of Theorem \ref{theo:construct}
we need some preparation.
First we describe a general construction of  
a knot $K_L^{\bar q}$
from a Seifert surface of $K'$, an $n$-component string link
$L$, and an $n$-tuple of integers ${\bar q}:=(q_1, \ldots, q_n)$.
Let $S' \subset \R ^{3}$ be a minimum genus  
Seifert surface for $K'$ and set
$g:={\rm genus}(S')$. Suppose that $S'$ is isotoped
into a disc-band form toward a spine $W_g$, which is a bouquet
of $2g$-circles based at a point $p$. Consider a  projection 
$P:\R ^{3} \longrightarrow R$ onto a projection plane $R$, 
so that the restrictions
of $P$ to $K'$ and $W_g$ are both regular.
We will identify $W_g$ with its diagram under the projection $P$. 
Let $D\subset R$ be a disc neighborhood of $p$, which contains no
crossing points of $W_g$.
Then, $D$ intersects  $W_g$ in a bouquet
of $4g$ arcs and the rest of $W_g$ consists of $2g$ arcs outside $D$.
We may assume that $S'$ is obtained from $W_g$ by replacing
each of the arcs outside $D$ by a band. 
Let $\alpha \subset \partial D$ be a connected subarc containing
$W_g\cap \partial D$ and set 
$\alpha':= \partial D\setminus \alpha$.
Let 
$$L: (I_n , \partial I_n) \longrightarrow 
(\overline{\R ^{3}\setminus D\times[0,1]}, \
\alpha')$$
be an $n$-component string link with components 
$L_1, \ldots, L_n$, where $I_n$ denotes the disjoint union
of $n$-copies of $I:=[0,1]$ and $D=D\times \{ \frac12 \}$.
The end points of $L_i$ in $\alpha'$ can be joined by a subarc $a_i$ in 
$\alpha'$ and we assume that $a_i\cap a_j=\emptyset$ if $i\neq j$.
Furthermore, 
we will assume that the restriction of $P$ on $L$ is regular
and the framing on $L$ defined by
parallel copies of $P(L)$ on $R$ is the zero
framing. If each $L_i$ is a subarc of $\alpha'$, we say that
the string link $L$ 
is {\it the trivial string link}. A string link $L$ in 
$\overline{\R ^{3}\setminus D\times[0,1]}$ is {\it trivial} 
if it is isotopic to
the trivial string link in $\overline{\R ^{3}\setminus D\times[0,1]}$ relative
to $\partial L$.  

We construct
a bouquet of $n+2g$ circles as follows:
For $1\leq i\leq n$, let $p_i, p_i'$ denote the endpoints of $P(L_i)$.
Connect
$p_i, p_i'$ to $p$ by disjointly embedded
arcs $\alpha_i, \alpha'_i$ that lie in $D$
and do not separate any of the arcs in
$D\cap W_g$. 
This process yields a bouquet $W_1:=W_1(L, W_g)$
of $n+2g$ circles.
Note that $W_1$ contains a sub-bouquet, say $W_L$, whose circles 
correspond to the components of $L$.

Let
${\bar q}:= (q_1, \ldots, q_n)$ be an $n$-tuple of integers.
For the circle in $W_L$ that corresponds to
the component $L_i$ of $L$ we add to $W_1$ an unlinked and unknotted
loop $L_{i}^{'}$,
which contains $q_i$ kinks. 
This is done in such
a way so that the four arcs of $L_i$ and 
$L_{i}^{'}$ in $D$
appear in alternating order. See
Figure 1. 
This produces a bouquet $W$ of $2(n+g)$ circles such that
$D\cap W$ is a bouquet
of $4(n+g)$ arcs and there are $2(n+g)$ arcs outside $D$.
Now we obtain a surface $S_L^{{\bar q}}$ by replacing
each of the arcs outside $D$ by a band, with twists
replacing the kinks contained on the arc. Let
$K_L^{{\bar q}}:=\partial S_L^{{\bar q}}$. If there is no danger 
of confusion we
will simply use $K_L$ to denote any of the knots $K_L^{{\bar q}}$. 
Next we prove two  lemmas needed for the proof
of Theorem \ref{theo:construct}.
%\vfill
%\eject
\begin{figure}
\centerline{\psfig{figure=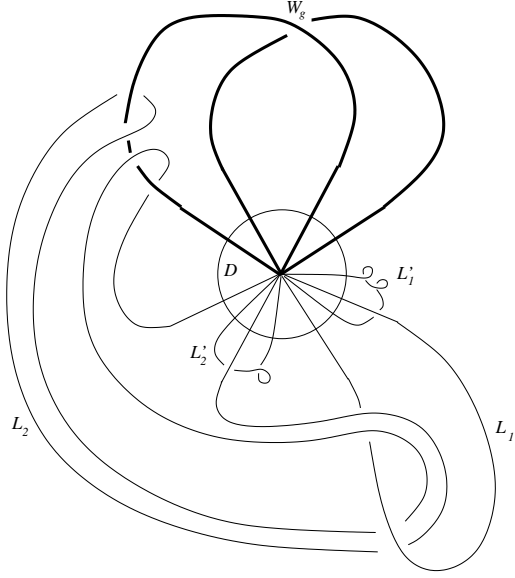,height=3in, clip=}}
\caption{An example of a bouquet $W$ and $D\cap W$}
\end{figure}

\begin{lem} \label{lem:sequivalent} Let $n>1$.
Suppose that the string link $L$ has the following
property: 
Every proper sublink $L'\subset L$ can be isotoped, relative to
$\partial L'$, in the complement of $W_g$ in 
$\overline{\R ^{3}\setminus D\times[0,1]}$, to the trivial string link
in $\overline{\R ^{3}\setminus D\times[0,1]}$.
Then, $K_L^{\bar q}$ is $n$-adjacent to $K'$ and it shares a common equivalence class
of Seifert matrices with $K'$.
Thus, in particular, $K_L^{\bar q}$ and $K'$ have the same Alexander module.
\end{lem}

\proof By construction of $S_L^{\bar q}$, for $1\leq i\leq n$,
the kinks on $L_i'$ give rise to a
generalized crossing, say $C_i$, of order $q_i$ on $K$.
Performing the generalized crossing changes in any non-empty subset of
$\{C_1,C_2,\dots,C_n\}$ will change $K_L^{\bar q}$ to $K_{L'}$ where $L'$ is
a proper subset of $L$. By the assumption on $L'$, $K_{L'}$ is isotopic to 
$K'$. Thus $K_L^{\bar q}$ is $n$-adjacent to $K'$.

Let $V$ denote the Seifert matrix of $S_L$ corresponding to the
spine $W$, and $V'$ denote the Seifert matrix of $S'$ corresponding to the
spine $W_g$. Since the linking number of $L_i$ with each circle in $W_g$ is 
zero, and the linking numbers between $L_i$ and $L_j$ are all zero, 
we see that 
$$V=\left[\begin{matrix}V' & 0 & 0   & \dots & 0 & 0 \\
                        0  & 0 & 0   & \dots & 0 & 0 \\
                        0  & 1 & q_1 & \dots & 0 & 0 \\
                         &\dots & \dots &\dots & \dots  & \\
                        0  & 0 & 0   &  \dots & 0 & 0 \\ 
                        0 & 0 & 0& \dots &      1 & q_n
\end{matrix}\right]$$
Thus, $K_L^{\bar q}$ and $K'$ have the same Alexander module. \qed
\smallskip
\smallskip
\smallskip

For the rest of the section we will assume that $L$ is chosen to satisfy 
the hypothesis
of Lemma \ref{lem:sequivalent}.
The next lemma describes the circumstances under 
which $K_L^{\bar q}$ is isotopic to $K'$.

\begin{lem} \label{lem:nonisotopic} Let $K'$ be a fibered knot.
Let $L$ be a string link in $\overline{\R^3\setminus D\times[0,1]}$
as in Lemma \ref{lem:sequivalent}.
Then, if $K_L^{\bar q}$ is isotopic to $K'$, $L$
can be isotoped, relative to
$\partial L$, in the complement of $W_g$ in 
$\overline{\R ^{3}\setminus D\times[0,1]}$, to the trivial string link
in $\overline{\R ^{3}\setminus D\times[0,1]}$.
\end{lem}

\proof
Suppose that $K_L^{\bar q}$ is isotopic to $K'$. For $i=1,\ldots, n$, let $D_i$
be a crossing disc corresponding to $C_i$ and
let $K_i:=\partial D_i$.
We can choose $K_i$ to be a small circle linking
once around the band of $S_L^{\bar q}$ corresponding to $L'_i$.
Let ${\bar K}_L^{\bar q}$ denote the knot
obtained from $K_L$ by changing all the generalized crossings
$C_1,\ldots, C_n$, simultaneously. One can see that
${\bar K}_L$ is obtained from $K'$
by $n$ finger moves, one for each component of $L$.
More specifically, to obtain ${\bar K}_L$,  for $i=1,\ldots, n$,
one pushes a small part of $K'=\partial S'$ that contains one endpoint of 
$L_i$ , following $L_i$ until one is getting very close to the other endpoint
of $L_i$. Then $K_L^{\bar q}$ is obtained
by allowing these fingers to intersect ${\bar K}_L$ so as to create
the generalized crossings $C_1,\ldots, C_n$.
See Figure 2.
It follows that  $D_1, \ldots, D_n$ are also crossing discs for
$K'$.

By Theorem 3.1 of \cite{kn:kl1}, a Seifert surface for
$K_L^{\bar q}$ that is of minimum genus in the complement of
$K_1 \cup \ldots
\cup K_n$ has to be a minimum genus surface for $K_L^{\bar q}$.
Since we assumed that ${\rm genus}(K_L^{\bar q})
= {\rm genus}(K')=g$, we conclude that
$K_L^{\bar q}$ bounds a Seifert surface
of genus $g$ in the complement of $K_1 \cup \ldots
\cup K_n$. Since ${\bar K}_L$ is obtained
from  $K_L^{\bar q}$ by twisting along
$D_1,\ldots, D_n$, the links
${\bar K}_L \cup K_1 \cup \ldots
\cup K_n$ and $K_L^{\bar q} \cup K_1 \cup \ldots
\cup K_n$ have homeomorphic complements.
We conclude that 
${\bar K}_L$ bounds a Seifert surface $\Sigma$ of genus $g$ in
the complement of $K_1 \cup \ldots
\cup K_n$. Since $\Sigma$ is incompressible, by isotopy of $\Sigma$
relative to $\partial\Sigma={\bar K}_L$, we can 
arrange so that ${\Sigma}\cap D_i$
is a single arc $b_i$ properly embedded in $\Sigma$. Each arc $b_i$ is 
a ``short'' subarc of $L_i$.

Clearly, performing the isotopy
from $K'$ to ${\bar K}_L$ described earlier backwards
isotopes the graph $b_1\cup  \ldots\cup  b_n\cup{\bar K}_L$ onto $L\cup K'$.
This isotopy brings $\Sigma$ to a minimal genus Seifert surface $\Sigma'$ of
$K'$. The string link $L$ lies on $\Sigma'$ as proper arcs. 
Since $K'$ is fibered, it admits a unique
minimum genus Seifert surface up to isotopy leaving $K'$
fixed point wise (see, for example, \cite{kn:bz}).
So, $\Sigma'$ and $S'$ are isotopic relative to $K'$. Since $L$ is disjoint
from $S'$, we may assume that during the isotopy from $\Sigma'$ to $S'$,
$L$ never touches $S'$ except for the last moment when $\Sigma'$ 
and $S'$ become identical. The isotopy from $L$ to its image on $S'$ are in
the complement of $S'$ and relative to $K'$.

\begin{figure}
\centerline{\psfig{figure=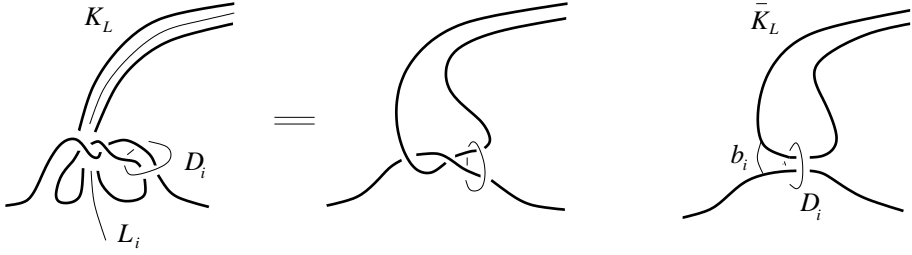,height=1.3in, clip=}}
\caption{ Constructing  $K_L^{\bar q}$ from ${\bar K}_L$}
\end{figure}

On the other hand, by  \cite{kn:k1},
each crossing $C_i$ must be nugatory.
If $C_i$ is a non-trivial nugatory crossing, we could have the
closure of the corresponding component of 
$L$ isotoped to a non-trivial summand in the connected sum decomposition
of $K'$. This contradicts to the assumption that each component of $L$ is 
trivial as a string link. Thus, we may assume that the image of $L$ on $S'$ 
is a collection of inessential proper arcs on $S'$. The position of 
the end points of this collection of inessential proper arcs force them to
bound disjoint disks on $S'$. 
We may assume that each of these disks lies in $D$, with one possible 
exception. For this exceptional disk, the corresponding proper arc in $S'$
would run out of $D$, follow the part of $K'$ outside of $D$, and come back 
to $D$. Then the closure of the corresponding component of $L$ would 
have been isotopic to $K'$. This is impossible. Thus the string link $L$ 
can be isotoped, relative to
$\partial L$, in the complement of $W_g$ in 
$\overline{\R ^{3}\setminus D\times[0,1]}$, to the trivial string link
in $\overline{\R ^{3}\setminus D\times[0,1]}$.
\qed

\vskip 0.08in

Next we turn our attention to the question of whether $K_L^{\bar q}$ can be chosen
to be prime. For this we need the following lemma:

\begin{lem} \label{lem:prime} Suppose that $K'$ is a prime knot and let $\bar q$ such that
$|q_i|>1$, for some $1\leq i\leq n$. Suppose, moreover, that
$g(K_L^{{\bar q}})>g(K')$, where $K_L^{{\bar q}}$ a knot associated to $L, {\bar q}$ as above.
If $K_L^{{\bar q}}$ is a composite knot then $K'$ is a summand of
$K_L^{{\bar q}}$ and  a 2-sphere realizing the corresponding decomposition of $K_L^{{\bar q}}$
separates $L$ from $W_g$.
\end{lem}

\proof
Suppose, without loss of generality, that $|q_1|>1$ and let
$K_1$ be a crossing link for the generalized crossing $C_1$ of 
$K_L^{{\bar q}}$.
Suppose that $K_L^{{\bar q}}$ has a non-trivial
connect sum composition $K_L^{{\bar q}}=J_1\# J_2$ and let $T$ be the corresponding
follow-swallow torus.
Since $g(K_L^{{\bar q}})>g(K')$, by Corollary 4.4 of \cite{kn:st1},  $K_1$ can be isotoped in the complement of $K_L^{{\bar q}}$ so that it is disjoint from $T$.
Let $V$ be the solid torus bounded by $T$; by assumption  $J_1$ is the core of $V$.
Suppose that $K_1$ lies outside $V$; the case that $K_1$ lies inside $V$ is completely analogous.
Then, $C_1$ is a crossing on $J_1$.
The knot obtained from $K_L^{{\bar q}}$ by changing  $C_1$
is of the form $J'_1\# J_2$,
where $J'_1$ the knot resulting from $J_1$.
By our assumptions on $L$,
$J'_1\# J_2=K'$. Since $K'$ is prime it follows that $J_2=K'$ and $J'_1$ is the unknot. Thus
$K_L^{{\bar q}}=J_1\# K'$. Let $Y$ be a 2-sphere that realizes this
connect sum. The surface $S'$ of $K'$
can be isotoped so that $S'\cap Y$
is an arc, say $\gamma$, properly embedded on $D$
such that the points in $\partial \gamma$
leave the arcs $\alpha, \alpha' \subset \partial D$
in different components of $S^3\setminus Y$.
Thus, in particular, $Y$ separates $\partial L$ from $W_g$.
Recall that $K_L^{{\bar q}}$ is the boundary of a surface obtained from
$S'\cup L$ by replacing each component of $L$ by an appropriate
band. Now,
since $K_L^{{\bar q}}\cap Y=\partial \gamma$, it follows that $Y$
separates $L$ from $W_g$. \qed
\smallskip
\vskip 0.04in

We can now finish the proof of Theorem \ref{theo:construct}:
\vskip 0.07in
\smallskip

{\bf Proof of Theorem \ref{theo:construct}:}
Let $n>1$. Let $S'$ be a minimum genus Seifert surface for $K'$
and a projection of it on $R$ as fixed earlier.  Choose 
$L: (I_n , \partial I_n) \longrightarrow ( \overline{ \R ^{3}\setminus D}, \alpha')$ such  that:

(i) Every proper sublink $L'\subset L$ can be isotoped, relatively $\partial L'$, in the complement of $W_g\cup K'$ so that it is properly embedded in $D$.

(ii) $L$ is not trivial.

(iii) There is no 2-sphere that intersects $K'$ at exactly two points and separates $L$ from $W_g$.

Let $\bar q$ such that
$|q_i|>1$, for some $1\leq i\leq n$. We claim that the knot $K:=K_L^{{\bar q}}$
has properties (a)-(d). Properties (a), (b), follow from (i) and Lemma \ref{lem:sequivalent}. By (ii) and Lemma \ref{lem:nonisotopic},
$K$ is not isotopic to $K'$; thus by Theorem \ref{theo:detect},
$g(K)>g(K')$. Now part (d) follows immediately from
Lemma \ref{lem:prime}. \qed

\smallskip
\appendix
\section{Obstructing symplectic structures}
In the recent years knots that look fibered to the Alexander
polynomial have received particular attention in symplectic geometry. For example,
a problem of current interest is when a 4-manifold of the form $S^1\times M$, where $M$ is a
3-manifold, admits a symplectic structure. It
is known that if $K$
is fibered then $S^1\times K(0)$ admits a
symplectic structure and it has been conjectured that the converse is true (see, \cite{kn:kr1} and references therein). It is known that the Alexander polynomial of a knot $K$ obstructs to the existence of
symplectic structures on $S^1\times K(0)$. More specifically, it is known that if
$S^1\times K(0)$ admits a
symplectic structure then $\Delta_K(t)$ is monic. Furthermore,
by a result of Kronheimer (\cite{kn:kr2}), if $g(K)>1$,  we must have $m(K)=0$.
Combining this with Proposition
\ref{pro:invariant}, it follows that  $\alpha(K)$ is a secondary obstruction to the
existence of
symplectic structures on $S^1\times K(0)$:

\begin{theo} \label{theo:symplectic} Let $K\in {\mathcal K}_1$. If
$S^1\times K(0)$ admits a symplectic structure
then, $\alpha(K)=0$.
\end{theo}
\proof If $g(K)>1$ the conclusion follows immediately from the aforementioned
result of Kronheimer and Proposition \ref{pro:invariant} (b).
Suppose that $g(K)=1$ and let $K'$ be a fibered knot such that
$\Delta_K(t)=\Delta_{K'}(t)$ and \ad , for some $n>1$.
By Theorem \ref{theo:detect},
$g(K)>g(K')$ and, since $g(K)=1$, $K'$ is the trivial knot. Hence
$\Delta_K(t)$ is trivial.
By Theorem 5.1 of \cite{kn:kl}, the only genus one knots that are
at least 2-adjacent to the unknot are
2-bridge knots. But the only 2-bridge knot with trivial Alexander polynomial is the unknot.
Thus $K=K'$, and by definition,
$\alpha(K)=0$. \qed
\smallskip

Let $K_n, K'$
be knots as in
Corollary \ref{corol:infinitely}. Since $K'$ is fibered,
$S^1\times K'(0)$ is symplectic. Since, by construction,
$\alpha(K_n)>0$, Theorem \ref{theo:symplectic} implies that
$S^1\times K_n(0)$ doesn't admit symplectic structures.
Thus we have examples of non-symplectic 4-manifolds
that are not distinguished from symplectic ones
by the information contained in
the Alexander polynomial.
Note that since $K_n(0), K'(0)$ have the same Alexander module,
it is know that
$S^1\times K_n(0)$ and
$S^1\times K'(0)$ are most distinguished by the
Seiberg-Witten invariants.


\smallskip






\begin{thebibliography}{XX10}
\smallskip

\bibitem[BZ]{kn:bz}
G. Burde and H. Zieschang, {\sl Knots.} de Gruyter Studies in Mathematics, 5. Walter de Gruyter \& Co., Berlin, 1985.


\bibitem[Ch]{kn:cha} Cha, Jae Choon, {\em Fibred knots and twisted Alexander invariants},
Trans. Amer. Math. Soc. 355 (2003), no. 10, 4187--4200.

\bibitem[C]{kn:co} T. Cochran, {\em Noncommutative knot theory},
 Algebr. Geom. Topol.  4  (2004), 347--398.

\bibitem[CM]{kn:cm} T. Cochran and P. Melvin, {\em Finite type invariants of 3-manifolds,} 
Invent. Math. 140 (2000), no. 1, 45--100.

\bibitem[CGLS]{kn:cgls}M. Culler, C. Gordon,
J. Luecke, and P. Shalen, {\em  Dehn surgery on knots,} Ann. of Math. (2) 125 (1987), no. 2, 237--300.

\bibitem[Ga1]{kn:ga} D. Gabai, {\em Foliations and the topology of 3-manifolds III},  J. Diff. Geom. vol {26} (1987), 479-536.


\bibitem[FKi]{kn:fk} S. Friedl and T. Kim, {\em Thurston norm, fibered manifolds and twisted Alexander polynomials,}
{\tt arXiv:math.GT/0505594}.


%\smallskip

\bibitem[Ga2]{kn:ga2} D.Gabai, {\em Detecting fibred links in $S\sp 3$}, Comment. Math. Helv. 61 (1986),
   no. 4, 519--555.

\bibitem[GKM]{kn:gkm} H.Goda, T.Kitano and T.Morifuji,
{\em Reidemeister torsion, twisted Alexander polynomial and fibered knots},
{\tt arXiv:math.GT/0311155.}





\bibitem[K]{kn:k} E. Kalfagianni, {\em Homology spheres with the same finite type invariants of bounded orders},  Math. Res. Lett. 4 (1997), no. 2-3, 341--347.


\bibitem[K1]{kn:k1} E. Kalfagianni, {\em Crossing changes of fibred knots,}
 preprint (2006), available at http://www.math.msu.edu/$\sim$kalfagia.


\bibitem[KL]{kn:kl} E. Kalfagianni and X.-S. Lin,
{\em Knot adjacency and satellites}, Topology and Applications,
 138 (2004), 207-217.

\bibitem[KL1]{kn:kl1} E. Kalfagianni and X.-S. Lin, {\em Knot adjacency,
genus and essential tori,} Pacific Jour. of Mathematics,  to appear. Posted at
{\tt arXiv:math.GT/0403024.}

\bibitem[Kr1]{kn:kr1} P. B. Kronheimer, {\em Embedded surfaces and gauge theory in three and four dimensions}, Surveys in differential geometry, Vol. III (Cambridge, MA, 1996), 243--298, Int. Press, Boston, MA, 1998.

\bibitem[Kr2]{kn:kr2} P. B. Kronheimer, {\em Minimal genus in $S\sp 1\times M\sp 3$}, Invent. Math. 135 (1999), no. 1, 45--61.


\bibitem[Mo]{kn:morton} H. R. Morton,
{\em Fibred knots with a given Alexander polynomial.}
{\sl Knots, braids and singularities} , 205--222,
Monogr. Enseign. Math., 31,
Enseignement Math., Geneva, 1983.




 \bibitem[O]{kn:o}T. Ohtsuki, {\em Finite type invariants of integral homology $3$-spheres,} J. Knot Theory Ramifications 5 (1996), no. 1, 101--115.


\bibitem[ST]{kn:st1} M. Scharlemann and A. Thompson, {\em
Unknotting number, genus, and companion tori,} Math. Ann. 280
(1988), no. 2, 191--205.



\bibitem[Tu]{kn:tu} V. Turaev, {\sl Torsions of $3$-dimensional manifolds.} Progress in Mathematics, 208. Birkhäuser Verlag, Basel, 2002.
%\smallskip
\end{thebibliography}
\end{document}